\documentclass[12pt,leqno]{article}
\usepackage{amsmath,amssymb,amsthm,amsfonts,mathrsfs,amscd,environ}
\usepackage{latexsym,enumerate,color,geometry}
\usepackage{fancyhdr,xcolor}
\usepackage{verbatim}
\usepackage{hyperref}
 \NewEnviron{ews}{%
\begin{equation}\begin{split}
  \BODY
\end{split}\end{equation}
}

\NewEnviron{ews*}{%
\begin{equation*}\begin{split}
  \BODY
\end{split}\end{equation*}
}

\def\beg{\begin}
\def\bequ{\begin{equation}}
\def\enqu{\end{equation}}
\def\bes{\begin{split}}
\def\ens{\end{split}}
\def\bews{\begin{ews}}
\def\beqn{\begin{eqnarray}}
\def\enqn{\end{eqnarray}}
\def\beq*{\begin{equation*}}
\def\enq*{\end{equation*}}
\def\bqn*{\begin{eqnarray*}}
\def\eqn*{\end{eqnarray*}}
\def\bary{\begin{array}}
\def\eary{\end{array}}
\def\bpma{\begin{pmatrix}}
\def\epma{\end{pmatrix}}
\def\bvma{\begin{Vmatrix}}
\def\evma{\end{Vmatrix}}

 \numberwithin{equation}{section}

\def\ga{\gamma}

\def\ep{\epsilon}

\def\th{\theta}

\def\la{\lambda}
\def\rh{\rho}

\def\si{\sigma}
\def\ta{\tau}
\def\ph{\phi}

\def\ps{\psi}

\def\Ga{\Gamma}

\def\De{\Delta}

\def\Ph{\Phi}
\def\Ps{\Psi}
\def\Om{\Omega}

\def\Q{\mathbb Q}
\def\R{\mathbb R}
\def\P{\mathbb P}
\def\E{\mathbb E}

\def\N{\mathbb N}

\def\sE{\mathscr E}
\def\sF{\mathscr F}
\def\sD{\mathscr D}

\def\sB{\mathscr B}
\def\sH{\mathscr H}

\def\cA{\mathcal A}
\def\cQ{\mathcal Q}
\def\cO{\mathcal O}

\def\d{\mathrm{d}}

\def\ff{\frac}
\def\ra{\rightarrow}
\def\nn{\nabla}
\def\pp{\partial}
\def\<{\langle}
\def\>{\rangle}
\def\sq{\sqrt}
\def\tld{\tilde}
\def\we{\wedge}

\def\Ran{\mathrm{Ran}}

\newcommand{\nov}[2]
{\begin{Vmatrix}\begin{pmatrix}
  #1\\#2
\end{pmatrix}\end{Vmatrix}}

\title{{\bf Harnack inequalities for $1$-$d$ stochastic Klein-Gordon type equations}\footnote{Supported by NSFC(11131003), SRFDP, 985-Project.}}

\author{
{\bf Shao-Qin Zhang }\\
\footnotesize{School of Statistics and Mathematics, Central University of Finance and Economics, Beijing 100081, China}\\
\footnotesize{Email: zhangsq@cufe.edu.cn}\\
\footnotesize{School of Math. Sci. and Lab. Math. Com. Sys., Beijing Normal University, Beijing 100875, China}\\
\footnotesize{Email: zhangsq@mail.bnu.edu.cn}
}

\begin{document}
\maketitle

\begin{abstract}
By the coupling method, we establish the Harnack inequalities, derivative formula and Driver's integration by parts formula for the stochastic Klein-Gordon type equations in the interval. We provide a detailed discussion about the nonlinear term. Some applications are given.
\end{abstract}\noindent

\noindent


\vskip 2cm

\section{Introduction}
The wave equation is the mathematical description of wave phenomena in physics. Nonlinear wave equations have been extensively study for many years, see \cite{Sog08} and reference thereby. When the wave motion is turbulent by random force, it is nature to consider the related model called stochastic
wave equations, see \cite{Chow02,DPZ1992,Dal09}. In this paper, we concern the following stochastic wave equation on an interval $\cO$ of $\R^1$:
\bequ\label{eq:NonlinearK-G0}
\beg{cases}
\d \dot X=\De X(t)\d t -l(X(t))\d t-\dot X(t)\d t+\si\d W(t),\\
X(0)=x\in H_0^1(\cO),\ \dot X(0)=y\in L^2(\cO),\\
X(t)=0, \mbox{~on~} \pp \cO,
\end{cases}
\enqu
where $\{W(t)\}_{t\geq 0}$ is a cylindrical Winer process on $L^2(\cO)$ in a complete filtered probability spaces $(\Om,\sF,\P, \{\sF_t\}_{t\geq 0})$. $\De$ is the Dirichlet-Laplace operator with $\sD(\De)=H^2(\cO)\cap H_0^1(\cO)$. $\dot X(t)=\ff {\d X(t)} {\d t}$ and $\si$ is a  Hilbert-Schmidt operator on $L^2(\cO)$. $l\in C(\R)$ satisfying the following conditions
\bequ
\beg{cases}
(1)\ l'\geq 0,\ |l(r)|\leq K_1|r|^\rh+K_2,\ |l'(r)|\leq K_3|r|^{\rh-1}+K_4;\\
(2)\ j(x):=\int_0^x l(r)\d r \geq K_5|x|^{\rh+1};\\
(3)\ |l'(r_1)-l'(r_2)|\leq \Big(C_1(|r_1|\we|r_2|)^{\rh-2}
+C_2\Big)|r_1-r_2|+C_3|r_1-r_2|^w,\ \rh> 2,\\
\hspace*{1.6em}|l'(r_1)-l'(r_2)|\leq C_4|r_1-r_2|^{\rh-1},\ \rh\in(1,2],\\
~~~~\ |l'(r_1)-l'(r_2)|\leq C_5(|r_1-r_2|^\ga\we1),\ \rh=1.
\end{cases}
\enqu
with $K_i(i=1,2,\cdots,5),\ C_i(i=1,\cdots,4)$ are some non-negative constants, $w\in(0,1)$, $\rh\geq 1$. Let $l(r)=|r|^{\rh-1}r$. Then $l$ satisfies (1)--(3) with $K_2=K_4=C_3=C_5=C_2=0$ and (\ref{eq:NonlinearK-G0}) is the stochastic Klein-Gordon Equation. Various problems had been concerned by many authors for stochastic wave equations. For example, \cite{Chow02,CarNua93,Crauel97} provided the existence and uniqueness of the solution of (\ref{eq:NonlinearK-G0}). \cite{CarNua93} concerned the existence of the random attractor. For ergodic properties, one can see \cite{BarbuD2002,DPZ1996}.

The dimension-free Harnack type inequalities was introduced by \cite{Wang97,Wang2013}. This type of inequalities have been established  not only  for various kinds of stochastic differential equations(SDEs) driven by Brownian motion (see \cite{WBook13',Wang2011,WangY2011,ShaoWY2012,GuiWang2012,Wang2013}), but also for SDEs driven by L\'evy noise and fractional Brownian motion, see \cite{WW,Wangjump2012,W13,OuyangRW2012,RoW2003,Fan13a,Fan13b}. Since the dimension-free property, the inequalities are possible valid for the infinite dimensional equations.  In fact, it has been proved that these inequalities holds for some stochastic partial differential equations(SPDEs), such as semilinear SPDEs, generalized porous media equations, fast-diffusion equations, stochastic Burgers equations and so on, see \cite{WBook13,LiuW2008,Liu2009,Ouyang2011,WangZTS2013,RoW2010,Zhang2013a,Zhang2013b} and reference there in. The derivative formula introduced in \cite{Bismut,ElwLiXM} and Driver's integration by part formula introduced in \cite{Driver} are both useful tools in stochastic analysis. They are closely linked to the Harnack type inequalities(see \cite{WBook13,W13}). The main aim of the paper is to establish  Harnack inequalities, the derivative formula and integration by part formula for the process $(X(t),\dot X(t))_{t\geq 0}$.

Though the stochastic wave equations can be rewrite as a semilinear SPDEs, we shall point out that it can not be covered by previous works. Let $H=L^2(\cO)$, $V=H_0^1(\cO)$. We denote $||\cdot||$ the norm of $L^2(\cO)$ and $||\cdot||_{H_0^1}$ the norm of $H_0^1(\cO)$. Then $\sH:=V\times H$ with norm $$\nov{x}{y}_\sH=\Big(||x||_{H_0^1}^2+||y||^2\Big)^{1/2}$$
is a Hilbert space. Let $Y(t)=\dot X(t)$ and
$$Z(t)=\bpma X(t)\\Y(t)\epma,\ \cA=\bpma 0&I\\ \De&0\epma,\ G(Z(t))=\bpma 0\\-l(X(t))-Y(t)\epma.$$
$\cA$ is an unbounded on $\sH$ with domain $\sD(\De)\times V$, moreover, it generates a $C_0$-group on $\sH$, saying
$$\Big(\bpma \cos(A^{1/2}t) & A^{-1/2}\sin(A^{1/2}t)\\-A^{1/2}\sin(A^{1/2}t) & \cos(A^{1/2}t)\epma\Big)_{t\geq 0}.$$
Then (\ref{eq:NonlinearK-G0}) can be write as the following semilinear SPDEs on $\sH$
\beq*
\beg{cases}
\d Z(t)=\cA Z(t)\d t+G(Z(t))\d t+\cQ\d W(t),\\
Z(0)=\bpma x\\y\epma,
\end{cases}
\enq*
where $\cQ$ is an operator from $H$ to $\sH$:
$$\cQ h=\bpma 0 \\\si h  \epma,\ \ \forall h\in H.$$
$\cQ$ is an injective map, but the range of $\cQ$(denoted by $\Ran(\cQ)$)~is~$\{0\}\times\Ran(\si)$, so $\bpma x\\y\epma\notin\Ran(\cQ)$. The Harnack inequalities for semilinear SPDEs established by previous works are usually dependent on the distance induced  by $(\cQ\cQ^*)^{-1/2}$(see \cite{RoW2003,RoW2010,WBook13,WangXu2012,WangZTS2013,Zhang2013a,Zhang2013b}). Therefore, we can not get the Harnack inequalities directly following the argument used previously. Recently, more and more works focus on degenerate SDEs and SPDEs, see \cite{WBook13'} and reference therein. \cite{GuiWang2012} introduced the coupling method to derive the Bismut formula for stochastic Hamilton systems. We extend the argument there to the stochastic wave equations and get the derivative formula. From the derivative formula and gradient-entropy inequality one can derive the Harnack inequality with power(see \cite{WBook13}), but in our situation, the derivative formula only holds for $\rh\in \{1\}\cup[2,\infty)$.  We start from the coupling again, and get the Harnack inequality with power for $\rh\in[1,2]$ just the same as the stochastic Hamilton system(see \cite{GuiWang2012}).

The paper is organized as follows. In Section 2, we first give some notation used frequently in the paper, and then state our main theorems and corollaries. We devote Section 3 to the proofs of our results.

\section{Main results}
We denote the Dirichlet-Laplace operator on $H$ by $-A$, then $A$ is a self adjoint operator on $H$. We endow the norm $||x||_{\th/2}:=||A^{\th/2}x||,\ x\in \sD(A^{\th/2})$ on the domain of $A^{\th/2}$. Then $||\cdot||_{H_0^1}=||\cdot||_{1/2}$.  Let~$P_Tf(x,y)=\E f(X_T(x), Y_T(y))$. We denote $\{e_j\}$ with the eigenvectors of $A$ and the eigenvalue of $A$ corresponding to $e_j$ by $\la_j$. Let $\si_0$ is a self-adjoint opertor on $H$ with $\si_0e_j=\si_{0j}e_j$ for some positive sequence $\{\si_0j\}$. Our first main result is
\beg{thm}\label{thm:derivative}
Assume that $\si\si^*\geq \si_0^2$, and there is $\la>0$ such that $\si_{0j}\sq{\la_j}\geq \ff 1 \la$. $\ga=1$ or $C_5=0$ if $\rh=1$, $C_3=0$ if $\rh>2$. Then, for all~$h_1\in \sD(A^{1/2}\si_0^{-1}),\ h_2\in \sD(\si_0^{-1})$, $\rh\in\{1\}\cup[2,\infty)$ and $v\in C^2([0,T],\R)$~with
$$v'(T)=v(T)=v'(0)=0,\ v(0)=1,$$
the derivative formula holds
\beg{ews*}
&\nn_{(h_1,h_2)} P_Tg(x,y)\\
&=\E g(X(T),Y(T))\int_0^T\Big\<\si^*(\si\si^*)^{-1}\Big[l'(X(t))\ps(t) +\ph(t)+f(t)\Big],\d W(t)\Big\>,\ g\in \sB(\sH),
\end{ews*}
where
\beg{ews*}
\ps(t) &=v(t)\Big(\cos(A^{1/2}t)h_1+A^{-1/2}\sin(A^{1/2}t)h_2\Big),\\
\ph(t)&= -v'(t)\Big(\cos(A^{1/2}t)h_1+A^{-1/2}\sin(A^{1/2}t)h_2)\\
&\qquad+v(t)\Big(\cos(A^{1/2}t)h_2-A^{1/2}\sin(A^{1/2}t)h_1\Big),\\
f(t)&=\bpma v''(t),\ 2v'(t)\epma\bpma \cos(A^{1/2}t) & A^{-1/2}\sin(A^{1/2}t)\\-A^{1/2}\sin(A^{1/2}t) & \cos(A^{1/2}t)\epma\bpma h_1\\ h_2\epma,\ t\in[0,T].
\end{ews*}
In fact, one can choose~$v(t)=1-\ff {3t^2} {T^2}+\ff {2t^3} {T^3}$.
\end{thm}

\beg{rem}
If $\si=A^{-1/2}$, then $\si^*(\si\si^*)^{-1}=A^{1/2}$ and $||\cdot||_{1/2}$ is equivalent to $||\cdot||_{H_0^1}$. Let $l(r)=|r|^{\rh-1}r$. If $\rh\in (1,2)$, then $l'(X(t))=|X(t)|^{\rh-1}$ is the fractional power of $|X(t)|$. So, it is not sure that $|X(t)|^{\rh-1}\in H_0^1(\cO)$ still holds. That means $\si^*(\si\si^*)^{-1}|X(t)|^{\rh-1}\ps(t)$ in the derivative formula may not make sense. Thus in Theorem \ref{thm:derivative}, $\rh\in \{1\}\cup[2,+\infty)$. It is slightly different from the finite dimensional case, see \cite[Theorem 2.2]{GuiWang2012}.
\end{rem}

Though we can not establish the derivative formula for $\rh\in (1,2)$, if starting from the coupling directly, it gives us a chance to obtain the Harnack inequality with power for $\rh\in [1,2]$ just as in finite dimensional case. If $\ker(\si\si^*)=\{0\}$, then $\si\si^*$ is positive operator. We can endow the space $\Ran(\si)$ with the norm $||x||_\si:=||(\si\si^*)^{-1/2}x||,\ x\in \Ran(\si)$. The norm on $L^p(\cO)$ is denoted by $||x||_{\rh+1}$ and
\beg{ews*}
\sE(x,y)&:=||x||_{H_0^1}^2+||y||^2+2J(x),\ J(x)=\int_{\cO}j(x(\xi))\d\xi.\\
\end{ews*}
Let $C_\cO$ be the Sobolev constant such that $\sup_{\cO}||\cdot||\leq C_\cO||\cdot||_{H_0^1}$ and 
\beg{ews*}
&\tld x=x+h_1,\ \tld y=y+h_2,\\ &|(h_1,h_2)|_{\si_0}=||\si_0^{-1}h_1||+||A^{-1/2}\si_0^{-1}h_2||,\\ &|(h_1,h_2)|_{\si_0+1/2}=||A^{1/2}\si_0^{-1}h_1||+||\si_0^{-1}h_2||,\\
&\sE_\si(p)=||\si||^2_{HS}+2(p-1)^+||\si||^2,\ \sE_T(p)=\ff {e^{(p-1)^+\sE_\si(p)T}-1} {(p-1)^+\sE_\si(p)T}.
\end{ews*}

\beg{thm}\label{thm:Harnack}
Assume that $\si\si^*\geq \si_0^2$, and there is $\la>0$ such that $\si_{0j}\sq{\la_j}\geq \ff 1 \la$.  Then, for all $h_1\in \sD(A^{1/2}\si_0^{-1}),\ h_2\in \sD(\si_0^{-1})$,\\
(1) for all~$\rh\geq 1$, $g\in \sB(\sH), g>0$, the~log-Harnack~inequality holds
\beg{ews*}
P_T \log g(\tld x,\tld y)\leq \log P_Tg(x,y)+\Ps_\rh(\tld x,\tld y,h_1, h_2,T\we1),
\end{ews*}
where
\beg{ews*}
&\Ps_\rh(\tld x,\tld y,h_1, h_2,T\we1)\\
&=\Ph_\rh(\tld x,\tld y,h_1,h_2,T\we1)+CK_1^2\Big[\ff {1+(T\we1)^2} {(T\we1)^3}|(h_1,h_2)|_{\si_0}^2+\ff {1+T\we1} {T\we1} |(h_1,h_2)|_{\si_0+1/2}^2\Big],
\end{ews*}
for $\rh=1$,
\beg{ews*}
\Ph_1(\tld x,\tld y,h_1,h_2,T\we1)&=\la^2(T\we1)\Big[(K_3+K_4)^2|(h_1,h_2)|_{1/2}^2\\
&\quad+C_5^2\Big(C_\cO^{2\ga}|(h_1,h_2)|_{1/2}^{2\ga}\we1\Big)\Big(|(h_1,h_2)|_{1/2}^2+\sE(\tld x,\tld y)+\sE_\si(1)(T\we 1)\Big),
\end{ews*}
for $\rh\in(1,2]$,
\beg{ews*}
&\Ph_\rh(\tld x,\tld y,h_1,h_2)\\
&=\la^2(T\we1)\Big\{|(h_1,h_2)|_{1/2}^2\Big[K_3^2C_{\cO}^{2\rh-2}\Big(\sE(\tld x, \tld y)+(\sE_\si(\rh-1)(T\we1))\Big)^{\rh-1}+K_4^2\Big]\\
&\quad+C_{\cO}^{2\rh-2}C_4^2\Big[|(h_1,h_2)|_{1/2}^{2\rh}+|(h_1,h_2)|_{1/2}^{2\rh-2}\Big(\sE( \tld x,\tld y)+\sE_\si(1)(T\we1)\Big)\Big]\Big\},
\end{ews*}
for $\rh\in(2,\infty)$,
\beg{ews*}
&\Ph_\rh(\tld x,\tld y,h_1,h_2)\\
&=\la^2C_{\cO}^{2\rh-2}(C_1^2+K_3^2)(T\we1)|(h_1,h_2)|_{1/2}^2\Big(\sE(\tld x,\tld y)+[\sE_\si(\rh-1)(T\we 1)]^{\ff {1} {\rh-1}}\Big)^{\rh-1}\sE_{T\we1}(\rh-1)\\
&\quad+2^{(\rh-1)\vee 2}\la^2C_{\cO}^{2\rh-2}(T\we1)\Big[C_3^2|(h_1,h_2)|_{1/2}^{2w+2}+C_1^2|(h_1,h_2)|_{1/2}^{2\rh-2}\Big(\sE( \tld x,\tld y)+\sE_\si(1)(T\we 1)\Big)\\
&\quad +C_1^2|(h_1,h_2)|_{1/2}^{2\rh}+C_2^2|(h_1,h_2)|_{1/2}^4\Big(\sE(\tld x, \tld y)+[\sE_\si(\rh-2)(T\we1)]^{\ff {1-(\rh-3)^-} {\rh-2}}\Big)^{\rh-2}\sE_{T\we1}(\rh-2)\Big]\\
&\quad+(T\we1)\Big[K_4^2|(h_1,h_2)|^2_{1/2}+C_3^2|(h_1,h_2)|^{2w}_{1/2}\Big(\sE(\tld x,\tld y)+\sE_{\si}(1)(T\we1)\Big)\Big].
\end{ews*}
with $C$~ an absolute constant.\\
(2) for all $\rh\in [1,2]$, the following Harnack inequality holds
$$(P_Tg(\tld x,\tld y))^p\leq P_Tg^p(x,y)\Ga(\tld x,\tld y,h_1,h_2),\ g\in \sB^+(H), p>1,$$
where for $\rh=1$, if $C_5=0$,
\beg{ews}
\Ga(\tld x,\tld y,h_1,h_2)= \exp\Big\{\ff {Cp} {(p-1)^2}\Big[\ff {1+T^2\we1} {T\we1}|(h_1,h_2)|_{1/2+\si_0}^2+\ff {1+T^2\we1} {T^3\we1} |(h_1,h_2)|_{\si_0}^2\Big]\Big\},
\end{ews}
and if $C_5>0$, then define $T_0=\ff {p-1} {4(C_5^2\vee1)\sq {2p} ||\si||}$,
\beg{ews}\label{eq:rh&a}
\Ga(\tld x,\tld y,h_1,h_2)&= \exp\Big\{\ff {Cp} {(p-1)^2}\Big[\ff {1+(T\we T_0)^2} {T\we T_0}|(h_1,h_2)|_{1/2+\si_0}^2+\ff {1+(T\we T_0)^2} {(T\we T_0)^3} |(h_1,h_2)|_{\si_0}^2\Big]\Big\}\\
&\quad\times\exp\Big\{(p-1)\Big(C_5^2C_\cO^{2\ga}|(h_1,h_2)|^{2\ga}_{1/2}\we1\Big)\Big[\ff {\sE(\tld x,\tld y)} {2||\si||^2(T\we T_0)} +\ff {||\si||_{HS}^2\log2} {||\si||^2}\Big]\Big\},
\end{ews}
for $\rh\in (1,2]$, let $T_0=\ff {\sq p-1} {4\sq 3||\si||\la C_\cO^{\rh-1}[\sq{K(h_1,h_2)}\vee1]}$,
\beg{ews*}
&\Ga(\tld x,\tld y,h_1,h_2)=\exp\Big\{(p-1)\ff {\Big((2-\rh)K_3^2+\ff {K_4^2} {C_\cO^{2\rh-2}}\Big)|(h_1,h_2)|_{1/2}^2
+C_4^2|(h_1,h_2)|_{1/2}^{2\rh}} {8||\si||^2(T\we T_0)[K(h_1,h_2)\vee1]}\Big\}\\
&\times\exp\Big\{\ff {Cp} {2(p-1)}\Big[\ff {1+(T\we T_0)^2} {T\we T_0}|(h_1,h_2)|_{1/2+\si_0}^2+\ff {1+(T\we T_0)^2} {(T\we T_0)^3}|(h_1,h_2)|_{\si_0}^2\Big]\Big\}\\
&\times \exp\Big\{\ff {2\sq p (\sq p+1)\tld c^2} {\sq p-1}[K(h_1,h_2)\we1]\Big(\ff {\sE(\tld x,\tld y)} {||\si||^2(T\we T_0)} +\ff {||\si||_{HS}^2\log4} {||\si||^2}\Big)\Big\},
\end{ews*}
where 
\bqn*
K(h_1,h_2)=K_3^2(\rh-1)|(h_1,h_2)|^2_{1/2}+C_4^2|(h_1,h_2)|_{1/2}^{2\rh-2}\\
\tld c^2=48||\si||^2(T\we T_0)^2\la^2C_\cO^{2\rh-2}[K(h_1,h_2)\vee1].
\eqn*
\end{thm}

Next, we shall consider Driver's integration by parts formula and shift-Harnack inequalitiess. Let $u\in C^2([0,T],\R)$, and
\beg{ews*}
\hat\ps(t) &=u(t)\Big(\cos(A^{1/2}(T-t))h_1+A^{-1/2}\sin(A^{1/2}(T-t))h_2\Big),\\
\hat\ph(t)&= u'(t)\Big(\cos(A^{1/2}(T-t))h_1+A^{-1/2}\sin(A^{1/2}(T-t))h_2)\\
&\qquad+u(t)\Big(A^{1/2}\sin(A^{1/2}(T-t))h_1-\cos(A^{1/2}(T-t))h_2\Big),\\
\hat f(t) &=\bpma u''(t),\ 2u'(t)\epma\bpma \cos(A^{1/2}(T-t)) & A^{-1/2}\sin(A^{1/2}(T-t))\\A^{1/2}\sin(A^{1/2}(T-t)) & -\cos(A^{1/2}(T-t))\epma\bpma h_1\\ h_2\epma.
\end{ews*}
Similar to Theorem \ref{thm:derivative} and Theorem \ref{thm:Harnack}, we have
\beg{cor}\label{cor:shift}
Assume that $\si\si^*\geq \si_0^2$, and there is $\la>0$ such that $\si_{0j}\sq{\la_j}\geq \ff 1 \la$. Let $h_1\in \sD(A^{1/2}\si_0^{-1}),\ h_2\in \sD(\si_0^{-1})$.\\
(1) Let $\ga=1$ or $C_5=0$ if $\rh=1$ and $C_3=0$ if $\rh>2$. Then for $\rh\in\{1\}\cup[2,\infty)$ and $u\in C^2([0,T],\R)$ satisfying
$$u(0)=u'(0)=u'(T)=0,\ u(T)=1,$$
the integration by parts formula holds
\beg{ews*}
&\nn_{(h_1,h_2)} P_Tg(x,y)\\
&=-\E g(X(T),Y(T))\int_0^T\Big\<\si^*(\si\si^*)^{-1}\Big[l'(X(t))\hat\ps(t) +\hat\ph(t)+\hat f(t)\Big],\d W(t)\Big\>,\ g\in C_b^1(\sH),
\end{ews*}
In fact, we can choose, for example, $u(t)=-\ff {3t^2} {T^2} + \ff {2t^3} {T^3}$.\\
(2) For $\rh\geq 1$, the shift log-Harnack inequality holds
$$P_T \log g(x,y)\leq \log P_T(g(\cdot+h_1,\cdot+h_2))(x,y)+\Ps_\rh(x,y,h_1, h_2,T),\ g>0,g\in\sB(\sH).$$
(3) For $\rh=1$ and $C_5=0$, the shift Harnack inequality holds
$$(P_Tg(x,y))^p\leq P_T(g^p(\cdot+h_1,\cdot+h_2))(x,y)\Ga(x,y,h_1,h_2),\ g\in \sB^+(\sH),$$
$\Ps_\rh,\ \Ga$ are the same as in Theorem \ref{thm:Harnack} except the  constant $C$.
\end{cor}
At last, we give some applications. It is almost standard, one can consult \cite{WBook13} for proofs.
\beg{cor}\label{cor:application}
Under the same assumption of Theorem \ref{thm:derivative}. Then\\
(1) for $\rh\in \{1\}\cup[2,+\infty)$ with $C_5=0$ if $\rh=1$ and $C_3=0$ if $\rh>2$, there exists $C>0$, such that the following gradient estimates hold
\beq*
|\nn P_Tg|^2(x,y)\leq C_1(1_{[\rh\geq 2]}\sE(x,y)+1)[P_Tg^2-(P_Tg)^2](x,y),\ (x,y)\in \sH.
\enq*
where $|\nn P_Tg|^2(x,y):=\sup_{||z||_{1/2+\si_0}\leq 1}\nn_z P_T g(x,y),$ and
$$\nn_z P_Tf(x,y):=\lim_{\ep\ra 0^+}\ff {P_T g((x,y)+\ep z)-P_T g(x,y)} {\ep}.$$
(2) let $P_T(z,\cdot),\ z\in \sH$ be the transition probability measure for $P_T$, then $P_T(z_1,\cdot)$ and $P_T(z_2,\cdot)$ are equivalent for $z_1,z_2\sH$ with $z_1-z_2\in \sD(A^{1/2}\si_0^{-1})\times\sD(\si_0^{-1})$. Let $p_{T,z_1,z_2}=\ff {\d P_T(z_1,\cdot)} {\d P_T(z_2,\cdot)}$. Then
\bqn*
&&P_T\{\log p_{T,z_1,z_2}\}(z_1)\leq \Psi_\rh(z_2,z_2-z_1,T\we1),\ \rh\in[1,+\infty),\\
&&P_T\{p^{\ff 1 {p-1}}_{T,z_1,z_2}\}(z_1)\leq \exp\Big(\ff {\Ga(z_2,z_2-z_1)} {p-1}\Big),\ \rh\in[1,2].
\eqn*
(3) $P_T(z,\cdot+z_1)$ is equivalent to $P_T(z,\cdot)$ for $z\in\sH$ and $z_1\in\sD(A^{(1+\th)/2})\times\sD(A^{\th/2})$. Let $p_{T,z,z_1}(y)=\ff {P_T(z,\d y+z_1)} {P_T(z,\d y )}$. Then
$$\int_{\sH}\exp\{p_{T,z,z_1}(y)\}P_T(z,\d y)\leq \exp\Big(\Psi_\rh(z+z_1,z_1)\Big).$$
\end{cor}

\section{The outline of proofs}
\emph{Proofs of Theorem \ref{thm:derivative} and Theorem \ref{thm:Harnack}}\\

Let $(\tld X(t),\tld Y(t))$ be the solution of equation, $\ep\in (0,1]$,
\bequ\label{eq:NonlinearK-Gc}\beg{cases}
\d \tld X(t)=\tld Y(t)\d t,\ \tld X(0)=x+\ep h_1,\\
\d \tld Y(t)=\Big(-A\tld X(t)-l(X(t))-Y(t)+\ep f(t)\Big)\d t\\
\hspace*{5em}+\si\d W(t),\ \tld Y(0)=y+\ep h_2.
\end{cases}
\enqu
Then, it is easy to see that
\bequ\beg{cases}\label{eq:X-tldX,Y-tldY}
\tld X(t)-X(t)=\ep \ps(t),\\
\tld Y(t)-Y(t)=\ep \ph(t),
\end{cases}\enqu
according to the definition of $\ps,\ph,f$. Particularly,
$$\tld X(T)=X(T),\ \tld Y(T)=Y(T).$$
Let
\beg{ews}\label{W_ep}
\d \tld W(t)&=\d W(t) + \si^*(\si\si^*)^{-1}\Big[l(\tld X(t))-l(X(t))+\ep \ph(t)+\ep f(t)\Big]\d t
\end{ews}
and
\beg{ews*}
R_s &=\exp\Big\{-\int_0^s\Big\<\si^*(\si\si^*)^{-1}\Big(l(\tld X(t))-l(X(t))+\ep \ph(t)+\ep f(t)\Big),\d W(t)\Big\>\\
&\hspace*{3em}-\ff 1 2 \int_0^s\Big|\Big|l(\tld X(t))-l(X(t))+\ep \ph(t)+\ep f(t)\Big|\Big|_\si^2\d t\Big\}.
\end{ews*}
Let $$\ta_n=\inf\{t\geq 0\ |\ ||X(t)||_{1/2}+||Y(t)||\geq n\}.$$
Then by Lemma \ref{Lmm:RlogR}, $\{R_{s\we\ta_n}\}_{n\in\N,s\in[0,T)}$ is uniformly integrable. By martingale convergence and domain convergence theorem, $\{R_s\}_{s\in[0,T]}$ is a martingale, moreover,
$$\E R_T\log R_T\leq \Ps_\rh(x,y,h_1, h_2).$$
According to the Girsanov theorem, $R_T\P$ is a probability measure and under $R_T\P$, $\{\tld W(t)\}_{t\in[0,T]}$ is cylindrical Brownian motion on $H$. So, $(\tld X(t),\tld Y(t))$ solves the following equation
\bequ\label{eq:NonlinearK-Gn}\beg{cases}
\d \tld X(t)=\tld Y(t)\d t,\ \tld X(0)=x+\ep h_1,\\
\d \tld Y(t)=-A\tld X(t)\d t-l(\tld X(t))\d t-\tld Y(t)\d t+\si\d \tld W(t),\ \tld Y(0)=y+\ep h_2,
\end{cases}
\enqu
Thus, $(\tld X(t),\tld Y(t))$ under the probability $R_T\P$ has the same distribution with $(X^{x+\ep h_1}(t),Y^{y+\ep h_2}(t))$ under $\P$, where $(X^{x+\ep h_1}(t),Y^{y+\ep h_2}(t))$ means solution of (\ref{eq:NonlinearK-G0}) with initial value $(x+\ep h_1,y+\ep h_2)$. By Lemma \ref{uintegrable}, and noting that $R_T$ is dependent on $\ep$,
\beg{ews*}
&\nn_{(h_1,h_2)}P_T g(x,y)=\lim_{\ep\ra 0^+}\ff {P_Tg(x+\ep h_1,y+\ep h_2)-P_Tg(x,y)} {\ep}\\
&=\lim_{\ep\ra 0^+} \ff  {\E R_Tg(\tld X(T),\tld Y(T))-\E g(X(T),Y(T))} {\ep}\\
&=\lim_{\ep\ra 0^+} \ff  {\E R_Tg( X(T), Y(T))-\E g(X(T),Y(T))} {\ep}\\
&=\E  g(X(T),Y(T))\ff {\d R_T} {\d \ep}|_{\ep=0}.
\end{ews*}
So, the derivative formula holds. Taking $\ep=1$, by the Young inequality, we obtain that
\beg{ews*}
&P_T\log g(x+h_1,y+h_2)=\E R_T\log g(\tld X(T),\tld Y(T))\\
&=\E R_T\log g(X(T),Y(T))\leq \log P_Tg(x,y)+\E R_T\log R_T,\ g\in \sB_b(\sH),\ g>0.
\end{ews*}
Letting $T_1\leq T$, due to the Markov property, we get that
\beg{ews*}
&P_T\log g(x+h_1,y+h_2)\leq P_{T_1}\log P_{T-T_1}g(x+h_1,y+h_2)\\
&\leq \log P_{T_1}P_{T-T_1}g(x,y) \E R_{T_1}\log R_{T_1}=\log P_Tg(x,y) \E R_{T_1}\log R_{T_1}.
\end{ews*}
Then part (1) of Theorem \ref{thm:Harnack} holds. Similarly, letting $T_1\leq T$, $\ep=1$
\beg{ews*}(P_Tf(\tld x,\tld y))^p&=(P_{T_1}(P_{T-T_1}f)(\tld x,\tld y))^p\leq P_{T_1}(P_{T-T_1}f)^p(x,y)(\E R^{p/{p-1}}_{T_1})^{p-1}\\
&\leq P_{T_1}P_{T-T_1}f^p(x,y)(\E R^{p/{p-1}}_{T_1})^{p-1}=P_Tf^p(x,y)(\E R^{p/{p-1}}_{T_1})^{p-1}.
\end{ews*}
So we only have to consider $T\leq \ff {\sq p-1} {4\sq 3||\si||c_0[(a\sq{K(h_1,h_2)})\vee1]}$. Therefore part (2) of Theorem \ref{thm:Harnack} follows from Lemma \ref{estimate:exp} and H\"older inequality.

The remainder of this section is devoted to the proofs of the technical lemmas. We start from a basic estimate of the energy $\sE(\tld X(t),\tld Y(t))$:
\beg{lem}\label{KL_lem_estimate}
For all $p\geq 1$, $s\in [0,T]$,
\beg{ews*}
\int_0^s\E R_{s\we\ta_n}\sE(\tld X(r\we\ta_n),\tld Y(r\we\ta_n))^p\d r\leq \Big(\sE(\tld x,\tld y)^p+\sE_\si(p)s\Big)\ff {e^{(p-1)\sE_\si(p)s}-1} {(p-1)\sE_\si(p)}
\end{ews*}
\end{lem}
\beg{proof} By~the It\^o~formula,
\beg{ews*}
\sE(\tld X(t),\tld Y(t))=\sE(\tld x,\tld y)-2\int_0^t||\tld Y(r)||^2\d r+2\int_0^t\<\tld Y(r),\si\d \tld W(r)\>+||\si||_{HS}^2t,\ t\leq s\we\tau_n.
\end{ews*}
So, for $p>1$
\beg{ews*}
&\d \sE(\tld X(t),\tld Y(t))^p\\
&\leq p\sE(\tld X(t),\tld Y(t))^{p-1}\Big(2\<\tld Y(t),\si\d \tld W(t)\>+||\si||_{HS}^2\d t\Big)\\
&\quad +2p(p-1)||\si||^2\sE(\tld X(t),\tld Y(t))^{p-1}\d t,\ t\leq s\we\tau_n.
\end{ews*}
Then by H\"older inequality
\beg{ews*}
&\E R_{s\we\ta_n}\sE(\tld X(t\we\ta_n),\tld Y(t\we\ta_n))^p\\
&\leq \sE(\tld x,\tld y)^p+\sE_\si(p)t+(p-1)\sE_\si(p)\int_0^t\E R_{s\we\ta_n}\sE(\tld X(r\we\ta_n),\tld Y(r\we\ta_n))^p\d r.
\end{ews*}
According to the Gronwall~lemma, we obtain that
\beg{ews*}
&\E R_{s\we\ta_n}\sE(\tld X(t\we\ta_n),\tld Y(t\we\ta_n))^p\\
&\leq \Big[\sE(\tld x,\tld y)^p+\sE_\si(p)t\Big]\exp{[(p-1)\sE_\si(p)t]},
\end{ews*}
and then the proof is  completed by integral from $0$ to $s$.
\end{proof}

\beg{lem}\label{Lmm:RlogR}
Under the same assumption of Theorem \ref{thm:Harnack}, then $$\sup_{s\in[0,T],n\geq 1}\E R_{s\we \ta_n}\log R_{s\we \ta_n} < \Ps_\rh(x,y,h_1, h_2).$$
\end{lem}
\beg{proof} By the definition of $R_s$,
\beg{ews*}
&\E R_{s\we \ta_n}\log R_{s\we \ta_n}\\
&\leq  \E R_{s\we \ta_n}\int_0^{s\we \ta_n}\Big|\Big|l(\tld X(t))-l(X(t))\Big|\Big|_\si^2\d t+\ep^2  \E R_{s\we \ta_n}\int_0^{s\we \ta_n}||\ph(t)+f(t)||_\si^2\d t
\end{ews*}
Since $\si\si^*\geq \si_0^2$, there exists an absolute constant $C$ such that
\beg{ews*}
&\E R_{s\we \ta_n}\int_0^{s\we \ta_n}||\ph(t)+f(t)||_\si^2\d t\\
&\leq C\Big[\ff {1+T^2\we1} {T^3\we1}|(h_1,h_2)|_{\si_0}^2+\ff {1+T^2\we1} {T\we1} |(h_1,h_2)|_{1/2+\si_0}^2\Big]
\end{ews*}
Since $||\si^{-1}_0 x||\leq \la ||x||_{H_0^1}$,
\beg{ews*}
&\E R_{s\we \ta_n}\int_0^{s\we \ta_n}\Big|\Big| l(\tld X(t))-l(X(t))\Big|\Big|_\si^2\d t\\
&\leq \la^2\E R_{s\we \ta_n}\int_0^{s\we \ta_n}\Big|\Big|\nn\Big(l(\tld X(t))-l(X(t))\Big)\Big|\Big|^2\d t\\
&=\la^2\E R_{s\we \ta_n}\int_0^{s\we \ta_n}\Big|\Big| l'(\tld X(t))(\nn\tld X(t)-\nn X(t))+(l'(\tld X(t))-l'(X(t)))\nn X(t) \Big|\Big|^2\d t.
\end{ews*}

If $\rh=1$, due to Sobolev's embedding theorem, $\sup_{\cO}|\cdot|\leq C_\cO ||\cdot||_{H_0^1}$, we have that
\beg{ews*}
&\E R_{s\we \ta_n}\int_0^{s\we \ta_n}\Big|\Big|l(\tld X(t))-l(X(t))\Big|\Big|_\si^2\d t\\
&\leq 3\la^2\E R_{s\we \ta_n}\int_0^{s\we \ta_n}\Big[(K_3+K_4)^2||\ps(t)||^2_{H_0^1}\\
&\qquad+C_5^2(\ep^{2\ga}C_\cO^{2\ga}||\ps(t)||_{H_0^1}^{2\ga}\we1)(||\ps(t)||_{H_0^1}^2
+||\tld X(t)||_{H_0^1}^2)\Big]\d t\\
&\leq 3\la^2(T\we1)\Big[(K_3+K_4)^2|(h_1,h_2)|_{1/2}^2\\
&\qquad+C_5^2\Big(\ep^{2\ga}C_\cO^{2\ga}|(h_1,h_2)|_{1/2}^{2\ga}\we1\Big)\Big(|(h_1,h_2)|_{1/2}^2+\sE(\tld x,\tld y)+\sE_\si(1)\Big)\Big]
\end{ews*}

If $\rh\in(1,2]$, we obtain that
\beg{ews*}
&\E R_{s\we \ta_n}\int_0^{s\we \ta_n}\Big|\Big|l(\tld X(t))-l(X(t))\Big|\Big|_\si^2\d t\\
&\leq \la^2\E R_{s\we \ta_n}\int_0^{s\we \ta_n}\Big|\Big|\ep (K_3|\tld X(t)|^{\rh-1}+K_4)|\nn\ps(t)|+\ep^\rh C_3|\ps(t)|^{\rh-1}|\nn\ps(t)|\\
&\hspace*{5em}+\ep^{\rh-1}C_3|\psi(t)|^{\rh-1}|\nn \tld X(t)|\Big|\Big|^2\d t\\
&\leq 3\la^2\E R_{s\we \ta_n}\int_0^{s\we \ta_n} \ep^2\Big(K_3C^{\rh-1}_{\cO}||\tld X(t)||^{\rh-1}_{H_0^1}+K_4\Big)^2||\ps(t)||_{H_0^1}^{2}\\
&\hspace*{1em}+C_4^2C_{\cO}^{2\rh-2}\Big(\ep^{2\rh}||\ps(t)||_{H_0^1}^{2\rh}+\ep^{2\rh-2}||\ps(t)||_{H_0^1}^{2\rh-2}||\tld X(t)||_{H_0^1}^{2}\Big)\d t.
\end{ews*}

If $\rh\in(2,\infty)$, then
\beg{ews*}
&\E R_{s\we \ta_n}\int_0^{s\we \ta_n}\Big|\Big|l(\tld X(t))-l(X(t))\Big|\Big|_\si^2\d t\\
&\leq \la^2\E R_{s\we \ta_n}\int_0^{s\we \ta_n}(K_3||\tld X(t)||_{H_0^1}^{\rh-1}+K4)^2||\ps(t)||_{H_0^1}^2\d t\\
&\quad+\la^2\E R_{s\we \ta_n}\int_0^{s\we \ta_n}\Big|\Big|\Big[\ep\Big(C_1(|\tld X(t)|\vee|X(t)|)^{\rh-2}+C_2\Big)|\ps(t)|\\
&\qquad+\ep^wC_3|\ps(t)|^w\Big]\Big[\ep|\nn\ps(t)|+
|\nn\tld X(t)|\Big]\Big|\Big|^2\d t.
\end{ews*}
According to $$(a\vee b)^p\leq 2^{(p-1)^+}((a\wedge b)^p+|b-a|^p),\  p>0,$$
and Sobolev's embedding theorem, we get that
\beg{ews*}
&\E R_{s\we \ta_n}\int_0^{s\we \ta_n}\Big|\Big|\ep^2(|\tld X(t)|\vee|X(t)|)^{\rh-2}|\ps(t)|\cdot|\nn\ps(t)|\Big|\Big|^2\d t\\
&\leq 2^{(\rh-3)^+}\E R_T\int_0^T\Big|\Big|\ep^2|\tld X(t)|^{\rh-2}|\ps(t)|\cdot|\nn\ps(t)|+\ep^{\rh-1}|\ps(t)|^{\rh-1}|\nn\ps(t)|\Big|\Big|^2\d t\\
&\leq 2^{(\rh-3)^++1}\E R_{s\we \ta_n}\int_0^{s\we \ta_n}\Big[\ep^4\sup_{\cO}\Big(|\tld X(t)|^{2\rh-4}|\ps(t)|^2\Big)||\nn\ps(t)||^2\\
&\hspace*{8em}+\ep^{2\rh}\sup_{\cO}\Big(|\ps(t)|^{2\rh-2}\Big)||\nn\ps(t)||^2\Big]\d t\\
&\leq 2^{(\rh-3)^++1}C_\cO^{2\rh-2}\E R_{s\we \ta_n}\int_0^{s\we \ta_n}\ep^{4}||\tld X(t)||_{H_0^1}^{2\rh-4}||\ps(t)||_{H_0^1}^4+\ep^{2\rh}||\ps(t)||_{H_0^1}^{2\rh}\d t.
\end{ews*}
Similarly,
\beg{ews*}
&\E R_{s\we \ta_n}\int_0^{s\we \ta_n}\Big|\Big|\ep(|\tld X(t)|\vee|X(t)|)^{\rh-2}|\ps(t)|\cdot|\nn\tld X(t)|\Big|\Big|^2\d t\\
&\leq 2^{(\rh-3)^++1}C_\cO^{2\rh-2}\E R_{s\we \ta_n}\int_0^{s\we \ta_n}\Big[\ep^2||\tld X(t)||_{H_0^1}^{2\rh-2}||\ps(t)||_{H_0^1}^2\\
&\hspace*{8em}+\ep^{2\rh-2}||\ps(t)||_{H_0^1}^{2\rh-2}||\nn\tld X(t)||^2\Big]\d t.
\end{ews*}
According to (\ref{eq:X-tldX,Y-tldY}), we have
\bequ
||\tld X(t)-X(t)||_{H_0^1(\cO)}\leq \ep|(h_1,h_2)|_{1/2}|v(t)|.
\enqu
Then, for $\rh>1$,
\beg{ews*}
&\E R_{s\we\ta_n}\int_0^{s\we\ta_n}||\tld X(t)||_{H_0^{1}}^{2\rh-2}v^2(t)\d t\\
&\leq \sup_{t\in[0,s]}v^2(t)\int_0^s\E R_{s\we\ta_n}||\tld X(t\we\ta_n)||_{H_0^1}^{2\rh-2}\d t\\
&\leq \sup_{t\in[0,s]}v^2(t)\Big(\sE(\tld x,\tld y)+(\sE_\si(\rh-1)s)^{\ff {1-(\rh-2)^-} {\rh-1}}\Big)^{\rh-1}\ff {e^{(\rh-2)^+\sE_\si(\rh-1)s}-1} {(\rh-2)^+\sE_\si(\rh-1)},
\end{ews*}
and
\beg{ews*}
&\E R_{s\we\ta_n}\int_0^{s\we\ta_n}||\tld X(t)||_{H_0^1}^2v^{2\rh-2}(t)\d t\\
&\leq \sup_{t\in[0,s]}v^{2\rh-2}(t)\int_0^{s}\E R_{s\we\ta_n}||\tld X(t\we\ta_n)||_{H_0^1}^2\d t\\
&\leq\sup_{t\in[0,s]}v^{2\rh-2}(t)\Big(\sE(\tld x, \tld y)+\sE_\si(1)s\Big)s.
\end{ews*}
For~$\rh>2$,
\beg{ews*}
&\E R_{s\we\ta_n}\int_0^{s\we\ta_n}||\tld X(t)||_{H_0^1}^{2\rh-4}v^4(t)\d t\\
&\leq \sup_{t\in[0,s]}v^4(t)\Big(\sE(\tld x,\tld y)+(\sE_\si(\rh-2)s)^{\ff {1-(\rh-3)^-} {\rh-2}}\Big)^{\rh-2}\ff {e^{(\rh-3)^+\sE_\si(\rh-2)s}-1} {(\rh-3)^+\sE_\si(\rh-2)}.
\end{ews*}
So, there exists a constant $C$ independent of $\ep$ such that
\beg{ews*}
&\E R_{s\we\ta_n}\log R_{s\we\ta_n}\leq \Ph_\rh(x,y,h_1,h_2)+C\Big[\ff {1+T^2\we1} {T^3\we1}|(h_1,h_2)|_{\si_0}^2+\ff {1+T^2\we1} {T\we1} |(h_1,h_2)|_{1/2+\si_0}^2\Big].
\end{ews*}
\end{proof}

\beg{lem}\label{uintegrable}
Under the assumptions of Theorem \ref{thm:derivative},
\beg{ews*}
\ff {\d} {\d \ep}|_{\ep=0}R^\ep_T=\int_0^T\Big\<\si^*(\si\si^*)^{-1}\Big[l'(X(t))\ps(t) +\ph(t)+f(t)\Big],\d W(t)\Big\>.
\end{ews*}
holds in $L^1(\P)$.
\end{lem}
\beg{proof} We follow~\cite{GuiWang2012}. It is clear that the lemma holds for $\rh=1$. So, we assume that $\rh\geq 2$. We write~$R^{(\ep)}_{T\we\ta_n},\ \tld X(t)^{(\ep)}$~and~$\tld W(t)^{(\ep)}$~to stress the parameter~$\ep$, see~(\ref{eq:X-tldX,Y-tldY})~and~(\ref{W_ep}). Since
\beg{ews*}
\ff {|R^{(\ep)}_{T\we\ta_n}-1|} \ep \leq &\ff 1 \ep \int_0^\ep R_{T\we\ta_n}^{(r)}\Big|\int_0^{T\we\ta_n}\<\si^*(\si\si^*)^{-1}\Big(l'(\tld X(t)^{(r)})\ps(t) +\ph(t)+f(t)\Big), \d W(t)\>\Big|\d r\\
&+\ff 1 \ep \int_0^\ep R^{(r)}_{T\we\ta_n}\Big(\int_0^{T\we\ta_n}\Big|\Big|l'(\tld X(t)^{(r)})\ps(t)+\ph(t)+f(t)\Big|\Big|_{\si}\\
& \times\Big|\Big|l(\tld X^{(r)}(t))-l(X^{(r)}(t))+r f(t)\Big|\Big|_{\si}\d t\Big)\d r,
\end{ews*}
we only have to prove that there exists a constant $C$ independent of~$n$~and~$r$ such that
\bqn*
&&\E R^r_{T\we\ta_n}\Big|\int_0^{T\we\ta_n}\<\si^*(\si\si^*)^{-1}\Big(l'(\tld X(t)^{(r)})\ps(t) +\ph(t)+f(t)\Big), \d W(t)\>\Big|^2<C\\
&&\E R^r_{T\we\ta_n}\int_0^{T\we\ta_n}\Big|\Big|l'(\tld X(t)^{(r)})\ps(t)+\ph(t)+f(t)\Big|\Big|^2_{\si}
d t<C\\
&&\E R^r_{T\we\ta_n}\int_0^{T\we\ta_n}\Big|\Big|l(\tld X^{(r)}(t))-l(X^{(r)}(t)) +r\ph(t)+r f(t)\Big|\Big|^2_{\si}\d t<C.
\eqn*
Noting that
\beg{ews*}
&\E R^r_{T\we\ta_n}\Big|\int_0^{T\we\ta_n}\<\si^*(\si\si^*)^{-1}\Big(l'(\tld X(t)^{(r)})\ps(t) +\ph(t)+f(t)\Big), \d W(t)\>\Big|^2\\
&\leq 2\E R^r_{T\we\ta_n}\Big|\int_0^{T\we\ta_n}\<\si^*(\si\si^*)^{-1}\Big(l'(\tld X(t)^{(r)})\ps(t) +\ph(t)+f(t)\Big), \d \tld W(t)^{(r)}\>\Big|^2\\
&+\E R^r_{T\we\ta_n}\int_0^{T\we\ta_n}\Big|\Big|l'(\tld X(t)^{(r)})\ps(t)+\ph(t)+f(t)\Big|\Big|^2_{\si}
d t\\
& +\E R^r_{T\we\ta_n}\int_0^{T\we\ta_n}\Big|\Big|l(\tld X^{(r)}(t))-l(X^{(r)}(t))  +r\ph(t)+r f(t)\Big|\Big|^2_{\si}\d t,
\end{ews*}
then we only have to prove
\bequ\label{K_G_inequ}\sup_{r\in(0,1),n\geq 1}\E R^r_{T\we\ta_n}\int_0^{T\we\ta_n}\Big|\Big|l'(\tld X^{(r)}(t))\ps(t)\Big|\Big|^2_{H_0^1}\d t<\infty.
\enqu
In fact, for $\rh=1$, since $\ga=1$ or $C_4=0$,
$$||l'(\tld X^{(r)}(t))||_{H_0^1}\leq C_4||\tld X^{(r)}(t)||_{H_0^1},$$
for~$\rh\geq 2$,
\beg{ews*}
\Big|\Big|l'(\tld X(t)^{(r)})\Big|\Big|_{H_0^1}&\leq \Big(C_1\sup_{\cO}|\tld X(t)^{(r)}|^{\rh-2}+C_2\Big)|\vee C_3|\tld X(t)^{(r)}||_{H_0^1}\\
&\leq C^{\rh-2}_\cO C_1||\tld X(t)^{(r)}||^{\rh-1}_{H_0^1}+C_2||\tld X(t)^{(r)}||_{H_0^1}.
\end{ews*}
Then there exists~$K'>0$, independent of $r,n$ such that
$$\Big|\Big|l'(\tld X(t)^{(r)})\ps(t)\Big|\Big|_{H_0^1}\leq K'||\ps(t)||_{H_0^1}||\tld X(t)^{(r)}||_{H_0^1}^{\rh-1}.$$
Therefore, (\ref{K_G_inequ}) holds due to Lemma \ref{KL_lem_estimate}.
\end{proof}

To obtain the estimate of $\E R_T^{p/(p-1)}$, we start form the following estimation, where we denote $\E_\Q$ the expectation w.r.t the probability measure $R_T\P$.
\beg{lem}
$$\E_\Q\exp\Big\{\ff 1 {8||\si||^2T^2}\int_0^T\sE(\tld X(t),\tld Y(t))\d t\Big\}\leq \exp\Big\{\ff {\sE(\tld x, \tld y)} {||\si||^2T}+\ff {||\si||^2_{HS}\log 4} {||\si||^2}\Big\}.$$
\end{lem}
\beg{proof}
Let
$$\ga(t)=\ff 1 {2||\si||^2(t+T)}.$$
Then
\beg{ews*}
\d (\ga(t)\sE(\tld X(t),\tld Y(t)))&=-2\ga(t)||\tld Y(t)||^2\d t+\ga'(t)\sE(\tld X(t),\tld Y(t))\d t\\
&\quad+2\ga(t)\<\tld Y(t),\si\d \tld W(t)\>+\ga(t)||\si||_{HS}^2\d t,
\end{ews*}
according to the It'\^o formula. So,
\beg{ews*}
&\E_\Q\exp\Big\{\int_0^s-\ga'(t)\sE(\tld X(t),\tld Y(t))\d t\Big\}\\
&\leq \exp\Big\{\ga(0)\sE(\tld x,\tld y)+||\si||_{HS}^2\int_0^s\ga(t)\d t\Big\}\E_\Q\exp\Big\{2\int_0^s\ga(t)\<\si^*\tld Y(t),\d \tld W(t)\>\Big\}\\
&\leq \exp\Big\{\ff {\sE(\tld x,\tld y)} {2||\si||_{HS}^2T} +\ff {||\si||_{HS}^2\log2 } {||\si||^2} \Big\}\Big(\E_\Q\exp\Big\{\int_0^s2\ga(t)^2||\si||^2\sE(\tld X(t),\tld Y(t))\d t\Big\}\Big)^{1/2}.
\end{ews*}
This implies that
$$\E_\Q\exp\Big\{\int_0^T\ff {\sE(\tld X(t),\tld Y(t))} {2||\si||^2(t+T)^2}\d t\Big\}\leq \exp\Big\{\ff {\sE(\tld x, \tld y)} {||\si||^2T} + \ff {||\si||_{HS}^2\log 4} {||\si||^2}\Big\}.$$
\end{proof}

Part (2) of Theorem \ref{thm:Harnack} follows from the lemma below.
\beg{lem}\label{estimate:exp}
Under the assumption of Theorem \ref{thm:Harnack} and moreover $T\leq \ff {\sq p-1} {4\sq 3||\si||\la C_\cO^{\rh-1}[\sq{K(h_1,h_2)}\vee1]}$ for  $\rh\in(1,2]$ and $T\leq \ff {p-1} {4(C_4^2\vee1)\sq {2p} ||\si||}$ for $\rh=1$, we have
$$(\E R_T^{\ff p {(p-1)}})^{p-1}\leq \Ga_\rh(\tld x,\tld y,h_1,h_2),\ h_1\in \sD(\si_0^{-1}A^{1/2}),\ h_2\in \sD(\si_0^{-1}).$$
\end{lem}
\beg{proof} First, we assume that $\rh\in (1,2]$. By the definition of $R_T$,
\beg{ews}
&\E R_T^{\ff p {p-1}}=\E R_T R_T^{\ff 1 {p-1}}=\E_\Q R_T^{\ff 1 {p-1}}\\
&\leq \E_\Q\exp\Big\{\ff {1} {(p-1)} \int_0^T\Big|\Big|l(\tld X(t))-l(X(t)) \Big|\Big|_\si^2\d t\\
&\quad+\ff 1 {p-1} \int_0^T\Big|\Big|\ph(t)+f(t)\Big|\Big|_\si^2\d t-\ff 1 {p-1} \int_0^T\<\si^*(\si\si^*)^{-1}\Big(\ph(t)+f(t)\Big),\d \tld W(t)\>\Big\}\\
&\quad-\ff {1} {p-1}\int_0^T\<\si^*(\si\si^*)^{-1}\Big(l(\tld X(t))-l(X(t)) \Big),\d \tld W(t)\>\\
&\leq \Big(\E_\Q\exp\Big\{\ff {2} {(p-1)}\int_0^T\Big|\Big|l(\tld X(t))-l(X(t)) \Big|\Big|_\si^2\d t\\
&\quad- \ff {2} {p-1} \int_0^T\<\si^*(\si\si^*)^{-1}\Big(l(\tld X(t))-l(X(t)) \Big),\d \tld W(t)\>\Big\}\Big)^{1/2}\\
&\quad\times \Big(\E_\Q\exp\Big\{\ff 2 {p-1} \int_0^T\Big|\Big|\ph(t)+f(t)\Big|\Big|_\si^2\d t\\
&\quad-\ff {2} {p-1}\int_0^T\<\si^*(\si\si^*)^{-1}\Big(\ph(t)+f(t)\Big),\d \tld W(t)\>\Big\}\Big)^{1/2}:=I_1^{1/2}\times I_2^{1/2}.
\end{ews}

Estimation of $I_2$:
\beg{ews*}
I_2&\leq \E_\Q\exp\Big\{\ff {-2} {p-1} \int_0^T\<\si^*(\si\si^*)^{-1}(\ph(t)+f(t)),\d \tld W(t)\>\\
&\quad -\ff 2 {(p-1)^2}\int_0^T||\ph(t)+f(t)||_\si^2\d t +\ff {2p} {(p-1)^2}\int_0^T||\ph(t)+f(t)||_\si^2\d t\Big\}\\
&\leq \exp\Big\{\ff {Cp} {(p-1)^2}\Big[\ff {1+T^2\we1} {T\we1}|(h_1,h_2)|_{1/2+\si_0}^2+\ff {(1+T^2\we1)} {T^3\we1} |(h_1,h_2)|_{\si_0}^2\Big]\Big\}.
\end{ews*}

For $I_1$. 
When $\rh\in (1,2]$,
\beg{ews*}
&\int_0^s\Big|\Big|l(\tld X(t))-l(X(t))\Big|\Big|_\si^2\d t\\
&\leq 3\la^2\int_0^s\Big(\Big(K_3^2C_\cO^{\rh-1}||\tld X(t)||_{H_0^1}^{\rh-1}+K_4)^2||\ps(t)||_{H_0^1}^2+C_4^2C_\cO^{2\rh-2}||\ps(t)||_{H_0^1}^{2\rh}\\
&\qquad+C_4^2C_\cO^{2\rh-2}||\tld X(t)||_{H_0^1}^{2}||\ps(t)||_{H_0^1}^{2\rh-2}\Big)\d t\\
&\leq 3\la^2\int_0^s\Big[(2-\rh)K_3^2C_\cO^{2\rh-2}||\ps(t)||_{H_0^1}^2
+C_4^2C_\cO^{2\rh-2}||\ps(t)||_{H_0^1}^{2\rh}+K_4^2||\ps(t)||_{H_0^2}^2\\
&\qquad + C_\cO^{2\rh-2}||\tld X(t)||_{H_0^1}^2\Big(K_3^2(\rh-1)||\ps(t)||_{H_0^1}^2+C_4^2||\ps(t)||_{H_0^1}^{2\rh-2}\Big)\Big]\d t.
\end{ews*}
Since $$K(h_1,h_2)=K_3^2(\rh-1)|(h_1,h_2)|^2_{1/2}+C_4^2|(h_1,h_2)|_{1/2}^{2\rh-2},$$
for $q>1$, we obtain that
\beg{ews}\label{ineq:I1}
I_1&\leq \Big(\E_\Q\exp\Big\{\ff {2(q-1+p)q} {(p-1)^2(q-1)} \int_0^s\Big|\Big|l(\tld X(t))-l(X(t))\Big|\Big|_\si^2\d t\Big\}\Big)^{\ff {q-1} {q}}\\
&\qquad\times\Big(\E_\Q\exp\Big\{-\ff {2q^2} {(p-1)^2}\int_0^s\Big|\Big|l(\tld X(t))-l(X(t))\Big|\Big|_\si^2\d t\\
&\qquad\qquad -\ff {2q} {p-1} \int_0^s\<\si^*(\si\si^*)^{-1}\Big(l(\tld X(t))-l(X(t))\Big),\d \tld W(t)\>\Big\}\Big)^{\ff 1 q}\\
&\leq \Big(\E_\Q\exp\Big\{\ff {2(q-1+p)q} {(p-1)^2(q-1)}\int_0^s\Big|\Big|l(\tld X(t))-l(X(t))\Big|\Big|_\si^2\d t\Big\}\Big)^{\ff {q-1} {q}}\\
&\leq \exp\Big\{\ff {6\la^2(q-1+p)} {(p-1)^2}\Big[(2-\rh)K_3^2C_\cO^{2\rh-2}\int_0^T||\ps(t)||_{H_0^1}^2\d t\\
&\qquad+C_4^2C_\cO^{2\rh-2}\int_0^T||\ps(t)||_{H_0^1}^{2\rh}\d t+K_4^2\int_0^T||\psi(t)||_{H_0^1}^2\d t\Big]\Big\}\\
&\qquad\times \Big(\E_\Q\exp\Big\{\ff {6\la^2C_\cO^{2\rh-2}(q-1+p)q} {(p-1)^2(q-1)}K(h_1,h_2)\int_0^T||\tld X(t)||_{H_0^1}^2\d t\Big\}\Big)^{\ff {q-1} {q}}
\end{ews}
Next, we shall estimate $$\E_\Q\exp\Big\{\ff {6\la^2C_\cO^{2\rh-2}(q-1+p)q} {(p-1)^2(q-1)}K(h_1,h_2)\int_0^T||\tld X(t)||_{H_0^1}^2\d t\Big\}.$$
Since
\bequ\label{ineq:exp1}
\ff {6\la^2C_\cO^{2\rh-2}(q-1+p)q} {(p-1)^2(q-1)}[K(h_1,h_2)\vee1]\leq \ff 1 {8||\si||^2T^2}
\enqu
is equivalent to
\bequ\label{eq:q}
(q-1)^2+\Big[(p+1)-\ff {(p-1)^2} {48||\si||^2T^2\la^2C_\cO^{2\rh-2}[K(h_1,h_2)\vee1]}\Big](q-1)+p\leq 0,
\enqu
for $T\leq \ff {\sq p-1} {4\sq 3||\si||\la C_\cO^{\rh-1}[\sq{K(h_1,h_2)}\vee1]}$, there exists $q>1$ such that (\ref{ineq:exp1}) holds. By Jensen's inequality,
\beg{ews}\label{ineq:exp3}
&\E_\Q\exp\Big\{\ff {6\la^2C_\cO^{2\rh-2}(q-1+p)q} {(p-1)^2(q-1)}K(h_1,h_2)\int_0^T||\tld X(t)||_{H_0^1}^2\d t\Big\}\\
&\leq \Big(\E_\Q\exp\Big\{\ff 1 {8||\si||^2T^2}\int_0^T||\tld X(t)||_{H_0^1}^2\d t\Big\}\Big)^{K(h_1,h_2)\we1}\\
&\leq \exp\Big\{[K(h_1,h_2)\we1]\Big(\ff {\sE(\tld x,\tld y)} {||\si||^2T} +\ff {||\si||_{HS}^2\log4} {||\si||^2}\Big)\Big\}.
\end{ews}

For other terms.
\beg{ews*}
&\exp\Big\{\ff {6\la^2(q-1+p)} {(p-1)^2}\Big[(2-\rh)K_3^2C_\cO^{2\rh-2}\int_0^T||\ps(t)||_{H_0^1}^2\d t\\
&\qquad+C_4^2C_\cO^{2\rh-2}\int_0^T||\ps(t)||_{H_0^1}^{2\rh}\d t+K_4^2\int_0^T||\psi(t)||_{H_0^1}^2\d t\Big]\Big\}\\
&\leq\exp\Big\{\ff {6(q-1+p)\la^2C_\cO^{2\rh-2}T} {(p-1)^2}\Big[\Big((2-\rh)K_3^2+\ff {K_4^2} {C_\cO^{2\rh-2}}\Big)|(h_1,h_2)|_{1/2}^2
+C_4^2|(h_1,h_2)|_{1/2}^{2\rh}\Big]\Big\}\\
&\leq \exp\Big\{\ff {\Big((2-\rh)K_3^2+\ff {K_4^2} {C_\cO^{2\rh-2}}\Big)|(h_1,h_2)|_{1/2}^2
+C_4^2|(h_1,h_2)|_{1/2}^{2\rh}} {8||\si||^2T[K(h_1,h_2)\vee1]}\Big\}.
\end{ews*}
Let $q$ be smallest the solution of (\ref{eq:q}) and $\tld c^2 = 48||\si||^2T^2\la^2C_\cO^{2\rh-2}[K(h_1,h_2)\vee1]$. Then
\beg{ews*}
&\ff {q-1} {q}(p-1)=\ff {4p(p-1)\Big(\ff {(p-1)^2} {\tld c^2} -(p+1)-\sq{\Big[\ff {(p-1)^2} {\tld c^2} -(p+1)\Big]^2-4p}\Big)} {\ff {(p-1)^2} {\tld c^2} -(p+1)-\sq{\Big[\ff {(p-1)^2} {\tld c^2} -(p+1)\Big]^2-4p}+1}\\
&\leq \ff {4p(p-1)\tld c^2} {(p-1)^2-\tld c^2(p+1)}\leq \ff {4p(p-1)\tld c^2} {2\sq p (\sq p-1)^2}=\ff {2\sq p (\sq p+1)\tld c^2} {\sq p-1},
\end{ews*}
where the two inequalities above we use that $\tld c^2\leq (\sq p-1)^2$. Combining this with (\ref{ineq:I1}), we get the estimation  for $\rh\in (1,2]$.
 
For $\rh=1$ and $C_5>0$. Since  $T\leq \ff {p-1} {4(C_5^2\vee1)\sq {2p} ||\si||}$, we have
$$\ff {4p(C_5^2\vee1)} {(p-1)^2} \leq \ff 1 {8||\si||^2T^2},$$
\beg{ews*}
I_1&\leq \Big\{\E_\Q\exp\Big\{\ff {4p} {(p-2)^2} \int_0^T||l(\tld X(t))-l(X(t))||_\si^2\d t\Big\}\Big\}^{1/2}\\
&\leq \exp\Big\{\ff {2pT} {(p-1)^2} \Big[(K_3+K_4)^2+C_5^2(C_\cO^{2\ga}|(h_1,h_2)|^{2\ga}\we1)|\Big](h_1,h_2)|_{1/2}^2\Big\}\\
&\qquad\times\Big\{\E_\Q\exp\Big\{\ff {4pC^2_5} {(p-1)^2}(C_\cO^{2\ga}|(h_1,h_2)|^{2\ga}_{1/2}\we1)\int_0^T||\tld X(t)||_{H_0^1}^2\d t\Big\}\Big\}^{1/2}\\
&\leq \exp\Big\{\Big(C_5^2C_\cO^{2\ga}|(h_1,h_2)|^{2\ga}_{1/2}\we1\Big)\Big[\ff {\sE(\tld x,\tld y)} {2||\si||^2T} +\ff {||\si||_{HS}^2\log2} {||\si||^2}\Big]\Big\}.
\end{ews*}
\end{proof}

\bigskip

\emph{Proof of Corollary \ref{cor:shift}}\\

We only have to consider the coupling for $(X(t),Y(t))_{t\geq0}$ and $(\hat X(t),\hat Y(t))_{t\geq0}$, where $(\hat X(t),\hat Y(t))_{t\geq0}$ is the solution of the following equation
\bequ\label{eq:NonlinearK-Gn-shift}\beg{cases}
\d \hat X(t)=\hat Y(t)\d t,\ \hat X(0)=x,\\
\d \hat Y(t)=-A\hat X(t)\d t-l(X(t))\d t\\
\hspace*{5em}-Y(t)\d t+\ep\hat f(t)\d t+\si\d  W(t),\ \hat Y(0)=y,
\end{cases}
\enqu
where $\ep\in(0,1]$. Repeating the argument in Theorem \ref{thm:derivative} and Theorem \ref{thm:Harnack}, one can get the corollary.
\qed

\bigskip

\textbf{Acknowledgement}
The author would like to thank Professor Feng-Yu Wang for useful suggestions.

\end{document}